\documentclass[11pt]{amsart}
\usepackage[margin=2.5cm]{geometry}
\usepackage{amsmath, amsthm, amssymb}
\usepackage{epsfig}
\usepackage{rawfonts}
\usepackage{enumerate}
\usepackage{graphics}
\usepackage{tikz}
\usepackage{tikz-qtree}
\usepackage{listings}
\usepackage{booktabs}

\usepackage{xspace}
     \usepackage{graphicx}

\theoremstyle{plain}

\theoremstyle{theorem}

\newtheorem{defn}{Definition}[section]
\newtheorem{prop}[defn]{Proposition}
\newtheorem{thm}[defn]{Theorem}
\newtheorem{lemma}[defn]{Lemma}
\newtheorem{conj}[defn]{Conjecture}

\newtheorem{coro}[defn]{Corollary}
\newtheorem{exa}[defn]{Example}
\newtheorem{rmk}[defn]{Remark}

\theoremstyle{remark}

\makeatletter
\@namedef{subjclassname@2020}{\textup{2020} Mathematics Subject Classification}
\makeatother

\begin{document}

\title[On some numerical semigroup transforms]{On some numerical semigroup transforms}
\author{Cisto Carmelo}

\address{Universit\'{a} di Messina, Dipartimento di Scienze Matematiche e Informatiche, Scienze Fisiche e Scienze della Terra\\
Viale Ferdinando Stagno D'Alcontres 31\\
98166 Messina, Italy}
\email{carmelo.cisto@unime.it}

		

\begin{abstract}

In this paper we introduce a particular semigroup transform $\mathcal{A}$ that fixes the invariants involved in Wilf's conjecture, except the embedding dimension. It also allows one to arrange the set of not ordinary and not irreducible numerical semigroups in a family of rooted trees.  We study also another transform, having similar features, that has been introduced by Bras-Amor\'os, and we make a comparison of them. In particular we study the behaviour of the embedding dimension under the action of such transforms, providing some consequences concerning Wilf's conjecture. 



\end{abstract}


\maketitle
\section{Introduction}

This paper is concerned with numerical semigroups. These are submonoids of $\mathbb{N}$ with finite complement in $\mathbb{N}$, a topic widely studied by several authors and from different perspectives. One of the nicest features that have been studied is the arrangement of the set of all numerical semigroups as a rooted tree, called the \emph{semigroup tree}. The semigroup tree has the set of all numerical semigroups as its set of vertices, and the set of all pairs $(S,S\cup \{\operatorname{F}(S)\})$ as its set of edges, where $\operatorname{F}(S)$ is the greatest element in $\mathbb{N}\setminus S$ (when $S\neq \mathbb{N}$). The root of this tree is $\mathbb{N}$. For these definitions see also \cite[Chapter 7, Section 1]{rosales2009numerical}. The building of the semigroup tree can be performed algorithmically and can be implemented in a programming language or by using a computer algebra software. We mention in particular the \texttt{GAP} \cite{GAP} package \texttt{numericalsgps} \cite{numericalsgps}, that contains many routines to deal with numerical semigroups. The semigroup tree is used in particular to produce all numerical semigroups of a given \emph{genus}, where the genus of a numerical semigroup $S$ is the number $|\mathbb{N}\setminus S|$. It is useful, for instance, to count the number $n_{g}$ of all numerical semigroups of a given genus $g$, as in \cite{bras2008fibonacci} which led to the
proposal of several interesting conjectures. Recall one of them, that states that the sequence $n_{g}$ has a Fibonacci-like behaviour, later proved in \cite{zhai2013fibonacci}. Another interesting argument related to numerical semigroups is the study of a conjecture posed for the first time by H. Wilf in \cite{wilf1978circle}. The semigroup tree has been used to test Wilf's conjecture up to a given genus, for instance up to genus 60 by Fromentin and Hivert in \cite{fromentin2016exploring}. Some improvements have been developed to explore the semigroup tree both from a computational point of view (for instance in \cite{bras2019right,fromentin2016exploring}) and from a theoretical point of view, depending on the properties one wants to examine (see \cite{bras2009bounds,bras2009towards} or the more recent \cite{delgado2019trimming}). The semigroup tree is related to the function on numerical semigroups defined by $S\mapsto S\cup \{\operatorname{F}(S)\}$. A function on numerical semigroups is called a \emph{numerical semigroup transform} (or simply a \emph{transform}). Given a transform it is possible, emulating the building of the semigroup tree, to arrange the set of numerical semigroups in a family of graphs, possibly trees. Thinking about Wilf's conjecture, we are interested in defining a transform that fixes two invariants involved in Wilf's conjecture and studying its consequences. The notion of \emph{special gap} and the technique that allows one to obtain the irreducible numerical semigroup that contains a fixed numerical semigroup (see \cite{rosales2003oversemigroups} and \cite[Chapter 3]{rosales2009numerical}) have inspired us to propose the definition of a transform with this particular features. We call such a transform $\mathcal{A}$ and we find certain numerical semigroups that will be the roots of the arranged rooted trees, and show they have a simple structure. This bring us to the definition of \emph{special} numerical semigroups. Moreover, Manuel Delgado who read a first preprint of this paper, pointed out to us that in \cite{Bras-Amoros2018inproc-Different} a transform having similar features had been introduced. We call such a transform $\mathcal{B}$ and we study some relationships between the two transforms.  \\
In Section 2 we summarize all notations, terminology and results needed to understand the rest of the paper. Section 3 is devoted to introducing special numerical semigroups and characterizing their structure. In Section 4 we define the transform $\mathcal{A}$, to which special numerical semigroups are related, providing some interesting properties of it, concerning in particular with the increasing of the embedding dimension. Subsequently we use the transform previously defined to arrange all non special numerical semigroups in a family of rooted trees, and this is the aim of Section 5 where we will also introduce some consequences related to Wilf's conjecture. In Section 6, we recall the transform $\mathcal{B}$ defined in \cite{Bras-Amoros2018inproc-Different} and we study it compared with $\mathcal{A}$. We conclude, in the last section, with some remarks and possible further developments.

\section{Preliminaries and known results}

Recall that a numerical semigroup $S$ is a submonoid of $\mathbb{N}$ such that $\mathbb{N}\setminus S$ is a finite set. It is well known that every numerical semigroup admits a unique finite minimal system of generators, that is, there exists a finite subset $\operatorname{G}(S)$ of $S$ such that every element of $S$ is obtained as a linear combination of elements in $\operatorname{G}(S)$ with coefficients in $\mathbb{N}$ and it is minimal in the sense that no proper subset of $\operatorname{G}(S)$ has the same property. The elements in $\operatorname{G}(S)$ are often called \emph{minimal generators}. Obviously an element $s\in S$ is not a minimal generator if and only if $s=s_{1}+s_{2}$ with $s_{1},s_{2}\in S\setminus \{0\}$. If a set $A$ generates a numerical semigroup $S$, in the sense described above, we usually write $S=\langle A\rangle$. Moreover every (minimal or not) system of generators  of a numerical semigroup is characterized by the fact that the greatest common divisor of all its elements is 1. For these and other interesting properties related to numerical semigroups a very good reference is \cite{rosales2009numerical}. If $S$ is a numerical semigroup, we provide
here several of its most important invariants that are useful for this paper:

\begin{itemize}
\item $\operatorname{H}(S)=\mathbb{N}\setminus S$ is called the set of \emph{gaps} of $S$.
\item $\operatorname{g}(S)=|\operatorname{H}(S)|$ is called the \emph{genus} of $S$.
\item $\operatorname{F}(S)=\max (\operatorname{H}(S))$ if $S\neq \mathbb{N}$, conventionally $\operatorname{F}(\mathbb{N})=-1$. It is called the \emph{Frobenius number} of $S$.
\item $\operatorname{m}(S)=\min (S\setminus \{0\})$ is called the \emph{muliplicity} of $S$.
\item $\operatorname{n}(S)=|\{s\in S\mid s<\operatorname{F}(S)\}|$, often referred to as the number of \emph{left elements} of $S$, if $S\neq \mathbb{N}$. Conventionally $\operatorname{n}(\mathbb{N})=1$.
\item $\operatorname{e}(S)=|\operatorname{G}(S)|$, the number of minimal generators, called the \emph{embedding dimension} of $S$. 
\end{itemize}

Observe that if for some $s \in S$ we have $\{s,s+1,\ldots,s+\operatorname{m}(S)-1\}\subset S$ then $s+n\in S$ for all $n\in \mathbb{N}$. We provide now some known results that we need for the forthcoming sections. The first one is quite easy to prove and we omit its proof.

\begin{prop}
Let $S$ be a numerical semigroup and $x\in S$. $S\setminus \{x\}$ is a numerical semigroup if and only if $x$ is a minimal generator of $S$.
\end{prop}

A numerical semigroup $S$ is called \emph{irreducible} if it cannot be expressed as an intersection of two numerical semigroups properly containing $S$. An irreducible numerical semigroup $S$ is called \emph{symmetric} if $\operatorname{F}(S)$ is odd, \emph{pseudo-symmetric} if $\operatorname{F}(S)$ is even. There are several characterizations for irreducible numerical semigroups. A useful result is the following:

\begin{prop}[\cite{rosales2009numerical}, Corollary 4.5]
Let $S$ be a numerical semigroup. Then
\begin{enumerate}
\item $S$ is symmetric if and only if $\operatorname{g}(S)=\frac{\operatorname{F}(S)+1}{2}$
\item $S$ is pseudo-symmetric if and only if $\operatorname{g}(S)=\frac{\operatorname{F}(S)+2}{2}$
\end{enumerate}

\end{prop}

\noindent Let $S$ be a numerical semigroup, we will use the following important subset of $\operatorname{H}(S)$:

$$\operatorname{SG}(S)=\{h\in \operatorname{H}(S)\mid 2h\in S,\ h+s\in S\ \mbox{for all}\ s\in S\}$$ 
$\operatorname{SG}(S)$ is called the set of \emph{special gaps} of $S$. The following nice results on special gaps can be found in \cite{rosales2003oversemigroups} or \cite{rosales2009numerical}.

\begin{prop}[\cite{rosales2003oversemigroups}, Proposition 9 and Corollary 13] Let $S$ be a numerical semigroup and $x\in \operatorname{H}(S)$. Then
\begin{enumerate}
\item $S\cup \{x\}$ is a numerical semigroup if and only if $x\in \operatorname{SG}(S)$.
\item $S$ is irreducible if and only if $\operatorname{SG}(S)=\{\operatorname{F}(S)\}$.
\end{enumerate}

\end{prop}

Some invariants of numerical semigroups are involved in a famous conjecture, widely studied by several authors: 

\begin{conj}[\textbf{Wilf's conjecture \cite{wilf1978circle}}] Let $S$ be a numerical semigroup. Then $$\operatorname{e}(S)\operatorname{n}(S)\geq \operatorname{F}(S)+1$$ or equivalently $$ (\operatorname{e}(S)-1)\operatorname{n}(S)\geq \operatorname{g}(S)$$ \label{Wilf}\end{conj}

It has been proved that Wilf's conjecture is satisfied by several classes of numerical semigroups, but it has not been proved to be true for every numerical semigroup. For a more complete and exhaustive survey about the study of Wilf's conjecture see~\cite{Delgado2020}.

\noindent Transforms on numerical semigroups is the main subject of this paper, in particular we are going to study two particular transforms. The definition of a numerical semigroup transform is expressed in the following: 

\begin{defn} \rm
Let $\mathcal{S}$ be the set of all numerical semigroups and $\mathcal{S}'\subseteq \mathcal{S}$. We call any function $\mathcal{F}:\mathcal{S}'\rightarrow \mathcal{S}$ a \emph{semigroup transform}. We denote $\mathcal{F}^{n}=\mathcal{F}^{n-1}\circ \mathcal{F}$, for $n>1$.
\end{defn}


\noindent We mention that there exist some particular transforms introduced in previous papers:
\begin{enumerate}

\item The \emph{``classical" transform} (\cite[Chapter 7]{rosales2009numerical}):
$$\mathcal{F}_{1}:\mathcal{S}\setminus \{\mathbb{N}\}\rightarrow \mathcal{S}\ \ \ \mbox{defined by}\ \ \ \mathcal{F}_{1}(S)=S\cup \{\operatorname{F}(S)\}.$$
In particular, for each $S\in \mathcal{S}\setminus \{\mathbb{N}\}$, there exists $n\in\mathbb{N}$ such that $\mathcal{F}_{1}^{n}(S)=\mathbb{N}$. This transform is related to the \emph{semigroup tree} (see for instance \cite{bras2009bounds,delgado2019trimming}).
 
\vspace{8pt} 
 
\item The \emph{ordinarization transform} (\cite{bras2012ordinarization}). A numerical semigroup $S$ is called $ordinary$ if there exists $c\in \mathbb{N}$ such that $S=\{s\in \mathbb{N}\mid s\geq c\}\cup \{0\}$, usually denoted by $S=\{0,c,\rightarrow\}$. Let $\mathcal{J}=\mathcal{S}\setminus \{S\in \mathcal{S}\mid S\ \mbox{is ordinary}\}$, we call $\mathcal{F}_{2}$ the transform:
$$\mathcal{F}_{2}:\mathcal{J}\rightarrow \mathcal{S}\ \ \ \mbox{defined by}\ \ \ \mathcal{F}_{2}(S)=(S\cup \{\operatorname{F}(S)\})\setminus \{\operatorname{m}(S)\}.$$
\noindent If $T=\mathcal{F}_{2}(S)$, then $\operatorname{g}(T)=\operatorname{g}(S)$. Moreover, for each $S\in \mathcal{S}$ there exists $n\in \mathbb{N}$ such that $\mathcal{F}_{2}^{n}(S)$ is the ordinary numerical semigroup of genus $\operatorname{g}(S)$. 

\vspace{8pt} 

\item The \emph{irreducibility transform} (\cite[Chapter 3]{rosales2009numerical}). Let $\mathcal{I}= \{S\in \mathcal{S}\mid S\ \mbox{is irreducible}\}$. For $S\in \mathcal{S}\setminus \mathcal{I}$ let $h=\max\{x\notin S\mid \operatorname{F}(S)-x\notin S, x\neq \frac{\operatorname{F}(S)}{2}\}$, we call $\mathcal{F}_{3}$ the transform:
$$\mathcal{F}_{3}:\mathcal{S}\setminus\mathcal{I}\rightarrow \mathcal{S}\ \ \ \mbox{defined by}\ \ \ \mathcal{F}_{3}(S)=S\cup \{h\}.$$

\noindent Observe that $S \subset \mathcal{F}_{3}(S)$ and that both semigroups have the same Frobenius number. Moreover, for each $S\in \mathcal{I}$ there exists $n\in \mathbb{N}$ such that $\mathcal{F}_{3}^{n}(S)$ is irreducible.

\end{enumerate}

\medskip

In the rest of this work if $a,b\in \mathbb{N}$ we denote $[a,b[\ =\{x\in \mathbb{N}\mid a\leq x<b\}$ and $[a,b]=\{x\in \mathbb{N}\mid a\leq x\leq b\}$.\\
For the sake of completeness we recall that an \emph{oriented graph} (or \emph{directed graph}) $G$ is a pair $(V,E)$, where $V$ is a nonempty set whose elements are called \emph{vertices}, and $E$ is a subset of $\{(v,w)\in V \times V\mid v\neq w \}$. The elements of $E$ are called \emph{edges} of $G$. A \emph{path} connecting the vertices $x$ and $y$ of $G$ is a sequence of distinct edges of the form $(v_0,v_1), (v_1,v_2),\ldots,(v_{n-1},v_n)$ with $v_0=x$ and $v_n= y$. An oriented graph $G$ is a \emph{rooted tree} if there exists a vertex $r$, known as the \emph{root} of $G$, such that for every other vertex $x$ of $G$, there exists a unique path connecting $x$ and $r$.

\section{Special numerical semigroups}

\begin{defn}\rm Let $S$ be a numerical semigroup. We say that $S$ is a $special$ numerical semigroup if $$(\operatorname{SG}(S)\setminus \{\operatorname{F}(S)\})\cap \{h\in \operatorname{H}(S)\mid h>\operatorname{m}(S)\}=\emptyset$$
\end{defn}

\noindent Standard examples of special numerical semigroups are the following:
\begin{itemize}
\item If $S$ is an irreducible numerical semigroup then $\operatorname{SG}(S)=\{\operatorname{F}(S)\}$, so it is special.
\item Recall that a numerical semigroup $S$ is called $ordinary$ if there exists $c\in \mathbb{N}$ such that $S=\{s\in \mathbb{N}\mid s\geq c\}\cup \{0\}$, denoted by $S=\{0,c,\rightarrow\}$. If $S$ is ordinary then $\{h\in \operatorname{H}(S)\mid h>\operatorname{m}(S)\}=\emptyset$, so it is special. Moreover $\operatorname{n}(S)=1$ and $\operatorname{e}(S)=c$.
\end{itemize}

Other special numerical semigroups are described in the next proposition.

\begin{prop}
Let $g\in \mathbb{N}$ with $g>2$ and $n\in [2,g]$. Then $S=\{0,g,g+1,\ldots,g+n-2,g+n,\rightarrow\}$ is a special numerical semigroup with $\operatorname{F}(S)=g+n-1$, $\operatorname{g}(S)=g$ and $\operatorname{n}(S)=n$. Furthermore:
\begin{enumerate}
\item If $n=g$ or $n=g-1$ then $S$ is irreducible.
\item If $3\leq n< g-1$ then $S$ is not irreducible and $\operatorname{e}(S)=g-1$.
\item If $g>3$ and $n=2$ then $S$ is not irreducible and $\operatorname{e}(S)=g$.
\end{enumerate}
\label{Sspecial}
\end{prop}

\begin{proof}
It is easy to verify that $S$ is a numerical semigroup with $\operatorname{F}(S)=g+n-1$, $\operatorname{g}(S)=g$ and $\operatorname{n}(S)=n$. Moreover it is special because for every $h\in \operatorname{H}(S)$ with $h\neq \operatorname{F}(S)$ we have $h<\operatorname{m}(S)=g$. \\
(1) If $n=g$ then $\operatorname{F}(S)+1=2g$, so $S$ is symmetric. If $n=g-1$ then $\operatorname{F}(S)+2=2g$, so $S$ is pseudo-symmetric (see \cite[Corollary 4.5]{rosales2009numerical}). In both cases it is irreducible. \\
(2) If $3\leq n<g-1$, $S$ is not irreducible, because $\operatorname{F}(S)-(g-1)=g+n-1-g+1=n\notin S$ (see \cite[Proposition 4.4]{rosales2009numerical}) and its minimal generators are $\{g,g+1,\ldots,g+n-2,g+n,g+n+1,\ldots,2g-1\}$, in fact $2g+j=g+(g+j)$ if $j\in \{0,1,\ldots,n-2,n\}$. \\
(3) If $n=2$ and $g>3$, $S$ is not irreducible because $\operatorname{F}(S)-(g-1)=2$ and it is easy to verify that its minimal generators are $\{g,g+2,\ldots,2g-1,2g+1\}$.
\end{proof}

In \cite{Bras-Amoros2018inproc-Different}, a numerical semigroup having only one gap greater then its multiplicity is called \emph{almost-ordinary}. Observe that such semigroups are exactly the numerical semigroups described in Proposition~\ref{Sspecial}. So, we use the same terminology and we highlight this fact with the following definition.

\begin{defn}\rm A numerical semigroup $S$ is called \emph{almost-ordinary} if $S=\{0,g,g+1,\ldots,g+n-2,g+n,\rightarrow\}$ such that $g\in \mathbb{N}$ with $g>2$ and $n\in [2,g]$.
\end{defn}

\begin{exa}\rm 
For $g=7$ and $n=4$, the numerical semigroup defined in Proposition~\ref{Sspecial} is the special numerical semigroup $S=\{0,7,8,9,11,\rightarrow\}=\langle 7,8,9,11,12,13\rangle$.
\end{exa}

For  $g=2$ the numerical semigroup $S$ in Proposition~\ref{Sspecial} is $S=\langle 2,5\rangle$, which is irreducible. Now we want to prove that every special numerical semigroup that is not irreducible nor ordinary is among those in Proposition~\ref{Sspecial} with $n<g-1$.

\begin{lemma}
Let $S$ be a numerical semigroup such that $\operatorname{F}(S)>2\operatorname{m}(S)$. Then one and only one of the following occurs:
\begin{enumerate}
\item $S$ is irreducible.
\item There exists $h\in \operatorname{SG}(S)\setminus \{\operatorname{F}(S)\}$ with $h>\operatorname{m}(S)$.
\end{enumerate}
\label{non-special}
\end{lemma}

\begin{proof}
If $S$ is irreducible then $\operatorname{SG}(S)=\{\operatorname{F}(S)\}$, so the second statement cannot occur. Suppose that $S$ is not irreducible, so there exist elements in $ \operatorname{SG}(S)\setminus \{\operatorname{F}(S)\}$ and it is easy to verify that such elements belong to $X=\left\lbrace x\notin S\mid \operatorname{F}(S)-x\notin S, x\neq \frac{\operatorname{F}(S)}{2}\right\rbrace$. In particular let $h=\max X$, then $h>\operatorname{F}(S)/2$. In fact if $h<\operatorname{F}(S)/2$, since $\operatorname{F}(S)-h\in X$ and $\operatorname{F}(S)-h>\operatorname{F}(S)/2>h$ we obtain a contradiction to the maximality of $h$ in $X$.\\ In particular $h>\operatorname{m}(S)$.
\end{proof}

\begin{thm}
Let $S$ be a special numerical semigroup. Suppose that $S$ is neither irreducible nor ordinary, then $\{s\in \mathbb{N}\mid \operatorname{m}(S)\leq s<\operatorname{F}(S)\}\subset S$.
\end{thm}

\begin{proof}
Let $k=\operatorname{F}(S)-\operatorname{m}(S)$, we prove that $\operatorname{F}(S)-j\in S$ for all $j=1,\ldots,k-1$. Suppose that $\operatorname{F}(S)-j\notin S$ for some $j\in\{1,\ldots,k-1\}$. We choose $j=\max\{i\mid 1\leq i\leq k-1, \operatorname{F}(S)-i\notin S\}$. $S$ is not ordinary, hence $\operatorname{m}(S)<\operatorname{F}(S)-j$. $S$ is special, so $\operatorname{F}(S)-j$ is not a special gap. It happens if one of the following occurs:
\begin{itemize}
\item[1)] $2(\operatorname{F}(S)-j)\notin S$ or
\item[2)] there exists $s\in S\setminus\{0\}$ such that $\operatorname{F}(S)-j+s\notin S$.
\end{itemize}
\noindent 1) In this case $2(\operatorname{F}(S)-j)=\operatorname{F}(S)$, that is, $\operatorname{F}(S)=2j$, so $\operatorname{F}(S)-(j-1)=j+1\in S$ and $\operatorname{m}(S)<\operatorname{F}(S)-j=j$. Let $n\in \mathbb{N}$ with $1\leq n<j$, such that $\operatorname{m}(S)=j-n$, then we have $\operatorname{F}(S)=2j=(j-n)+(j+n)=\operatorname{m}(S)+(j+n)$, that is, $j+n\notin S$. Moreover $\operatorname{F}(S)>j+n\geq j+1=\operatorname{F}(S)-(j-1)>\operatorname{F}(S)-j$, so $j+n\in S$, a contradiction.\\
2) In this case there exists $s\in S\setminus \{0\}$ such that $\operatorname{F}(S)-j+s=\operatorname{F}(S)$. We prove that it is a contradiction by showing that $\operatorname{F}(S)-j+\operatorname{m}(S)>\operatorname{F}(S)$ and considering that $s\geq \operatorname{m}(S)$. If $\operatorname{F}(S)-j+\operatorname{m}(S)=\operatorname{F}(S)$ then $\operatorname{m}(S)=j$ and $\operatorname{F}(S)>2\operatorname{m}(S)$, since $\operatorname{F}(S)>\operatorname{m}(S) +j$. So by Lemma~\ref{non-special}, $S$ is not special, contradicting our hypothesis. If $\operatorname{F}(S)-j+\operatorname{m}(S)<\operatorname{F}(S)$ then $\operatorname{m}(S)<j$, in particular $1<\operatorname{m}(S)\leq j-1$, so $\operatorname{F}(S)-\operatorname{m}(S)>\operatorname{F}(S)-j$. This means, by the choice of $j$, that $\operatorname{F}(S)-\operatorname{m}(S)\in S$, which is a contradiction. 
\end{proof}

\begin{coro}
Let $S$ be a special numerical semigroup with genus $g>3$. If $S$ is neither irreducible nor ordinary then $S=\{0,g,g+1,g+2,\ldots,g+n-2,g+n,\rightarrow\}$ for some $n \in \mathbb{N}$ with $2\leq n<g-1$, in particular $S$ is an \emph{almost-ordinary} numerical semigroup.
\end{coro}

\begin{proof}
Easily follows from the previous theorem and Proposition \ref{Sspecial}.
\end{proof}

\begin{prop}
Every special numerical semigroup satisfies Wilf's conjecture.
\end{prop}
\begin{proof}
Let $S$ be a special numerical semigroup. If $S$ is ordinary the claim easily follows. If $S$ is irreducible then $S$ satisfies Wilf's conjecture by \cite[Proposition 2.2]{dobbs2003question}. If $S$ is not irreducible then, by Proposition~\ref{Sspecial}, $2\operatorname{e}(S)\geq\operatorname{m}(S)$ so $S$ satisfies Wilf'conjecture by \cite[Theorem 18]{sammartano2012numerical}.
\end{proof}

\section{The transform $\mathcal{A}$}

\begin{defn}\rm Let $\mathcal{H}$ be the set of all non special numerical semigroups. We introduce the following transform:

$$\mathcal{A}:\mathcal{H}\rightarrow \mathcal{S}\ \ \ \mbox{defined by}\ \ \ \mathcal{A}(S)=(S\cup \{h\})\setminus \{\operatorname{m}(S)\}\ \ \ \mbox{where}\ \ \ h=\max (\operatorname{SG}(S)\setminus\{\operatorname{F}(S)\}).$$

\end{defn}

\noindent The transform $\mathcal{A}$ is actually a mixture of $\mathcal{F}_{2}$ and $\mathcal{F}_{3}$, considering also the following:

\begin{lemma}
Let $S$ be a non irreducible numerical semigroup. Then 
$$\max (\operatorname{SG}(S)\setminus \{\operatorname{F}(S)\})=\max\ \left\lbrace x\notin S\mid \operatorname{F}(S)-x\notin S, x\neq \frac{\operatorname{F}(S)}{2}\right\rbrace$$
and it is greater than $\frac{\operatorname{F}(S)}{2}$.
\label{max=max}
\end{lemma}

\begin{proof}
Let $h_{1}=\max (\operatorname{SG}(S)\setminus \{\operatorname{F}(S)\})$ and $h_{2}=\max \{ x\notin S\mid \operatorname{F}(S)-x\notin S, x\neq \frac{\operatorname{F}(S)}{2}\}$. From \cite[Lemma 4.3]{rosales2009numerical}, $h_{2}\in \operatorname{SG}(S)$ and $h_{2}>\frac{\operatorname{F}(S)}{2}$, so $h_{2}\leq h_{1}$. Suppose that $h_{2}<h_{1}$, then $h_{1}\notin \{ x\notin S\mid \operatorname{F}(S)-x\notin S, x\neq \frac{\operatorname{F}(S)}{2}\}$ and since $h_{1}\notin S$ we have $2h_{1}=\frac{\operatorname{F}(S)}{2}$ or $\operatorname{F}(S)-h_{1}\in S$, but this is a contradiction since $h_{1}$ is a special gap. So $h_{2}=h_{1
}$.
\end{proof}

\noindent We want to study some properties of the transform $\mathcal{A}$.

\begin{prop}
Let $S$ be a non special numerical semigroup and $T=\mathcal{A}(S)$. Then $\operatorname{F}(T)=\operatorname{F}(S)$, $\operatorname{g}(T)=\operatorname{g}(S)$ and $\operatorname{n}(S)=\operatorname{n}(T)$.

\end{prop}
\begin{proof}
The assertions easily follow from the definition of $\mathcal{A}$.
\end{proof}

\begin{lemma}
Let $S$ be a non special numerical semigroup and $T=\mathcal{A}(S)$. Then $T$ is neither irreducible nor ordinary.
\label{lemmanor}
\end{lemma}

\begin{proof}
Let $S$ be a non special semigroup. Then there exists $h=\max (\operatorname{SG}(S)\setminus \{\operatorname{F}(S)\})$ and $h>\operatorname{m}(S)$. We have $T=(S\cup \{h\})\setminus \{\operatorname{m}(S)\}$, in particular $T\cup \{\operatorname{m}(S)\}=S\cup \{h\}$ is a numerical semigroup, $\operatorname{m}(S)\neq \operatorname{F}(T)$ (since $\operatorname{F}(T)=\operatorname{F}(S)$), that is $\operatorname{m}(S)\in \operatorname{SG}(T)\setminus \{\operatorname{F}(T)\}$, so $T$ is not irreducible. $T$ is not ordinary because $h\in T$ and $h<\operatorname{F}(T)$.
\end{proof}

\begin{thm}
Let $S$ be a non special numerical semigroup. Then there exists $r\in \mathbb{N}$ such that $\mathcal{A}^{r}(S)=\{0,g,g+1,\ldots,g+n-2,g+n,\rightarrow\}$ with $n=\operatorname{n}(S)$ and $g=\operatorname{g}(S)$, that is, $\mathcal{A}^{r}(S)$ is an almost-ordinary numerical semigroup.
\label{rootNS}
\end{thm}

\begin{proof}
Let $S$ be a non special numerical semigroup. If $T=\mathcal{A}(S)$ then we have $\operatorname{m}(T)>\operatorname{m}(S)$, $\operatorname{F}(S)=\operatorname{F}(T)$ and $\operatorname{n}(S)=\operatorname{n}(T)$. In particular $|\{h\in \operatorname{H}(T)\mid h>\operatorname{m}(T)\}|< |\{h\in \operatorname{H}(S)\mid h>\operatorname{m}(S)\}|$. Therefore, if we apply repeatedly the transform $\mathcal{A}$, starting from $S$, after $r>0$ steps we obtain a special numerical semigroup, that is $\mathcal{A}^{r}(S)$ and, by Lemma~\ref{lemmanor}, it is almost-ordinary. 
\end{proof}

\begin{exa} \rm
Let $S=\{0,5,7,10,12,\rightarrow\}=\mathbb{N}\setminus \{1,2,3,4,6,8,9,11\}$.\\
We have $\max (\operatorname{SG}(S)\setminus\{\operatorname{F}(S)\})=9$ and $\operatorname{m}(S)=5$ so $\mathcal{A}(S)=\{0,7,9,10,12,\rightarrow\}=\mathbb{N}\setminus \{1,2,3,4,5,6,8,11\}$. \\
Furthermore $\mathcal{A}^{2}(S)=\{0,8,9,10,12,\rightarrow\}=\mathbb{N}\setminus \{1,2,3,4,5,6,7,11\}$ which is almost-ordinary.

\end{exa}

The transform $\mathcal{A}$ fixes all invariants involved in Wilf's conjecture except for the embedding dimension $\operatorname{e}(S)$. If $T=\mathcal{A}(S)$, we are interested in studying when $\operatorname{e}(T)$ increases or decreases with respect to $\operatorname{e}(S)$.

\begin{prop}
Let $S$ be a non special numerical semigroup, $h=\max (\operatorname{SG}(S)\setminus \{\operatorname{F}(S)\})$ and $T=\mathcal{A}(S)$. Then $2\operatorname{m}(S)$ and $h$ are minimal generators of $T$.
\label{h,2n1}
\end{prop}
\begin{proof}
It is evident that $h$ is a minimal generator of $T$. If $2\operatorname{m}(S)$ is not a minimal generator of $T$, then $2\operatorname{m}(S)=t_{1}+t_{2}$ with $t_{1},t_{2}\in T$. But $\operatorname{m}(T)>\operatorname{m}(S)$ so for all $t\in T$ we have $t>\operatorname{m}(S)$ and $2\operatorname{m}(S)=t_{1}+t_{2}>2\operatorname{m}(S)$, which is a contradiction. 
\end{proof}

\begin{lemma}
Let $S$ be a non special numerical semigroup, $h=\max (\operatorname{SG}(S)\setminus \{\operatorname{F}(S)\})$ and $T=\mathcal{A}(S)$. Let $\{\operatorname{m}(S)=n_{1}<n_{2}<\ldots<n_{t}\}$ be the set of minimal generators of $S$. If $\operatorname{m}(S)+h\geq n_{r}$ for some $r\in \{1,\ldots,t\}$, then $n_{2},n_{3},\ldots,n_{r}$ are minimal generators of $T$. In particular if $\operatorname{m}(S)+h\geq n_{t-1}$ then $\operatorname{e}(T)\geq \operatorname{e}(S)$.
\label{<nr}
\end{lemma}
\begin{proof}
Let $n_{i}\leq \operatorname{m}(S)+h$ and suppose that $n_{i}$ is not a minimal generator of $T$. Then $n_{i}=t_{1}+t_{2}$ with $t_{1},t_{2}\in T$. Since $n_{i}$ is a minimal generator of $S$ the only possibilities are that $n_{i}=s+h$ with $s\in S\setminus \{\operatorname{m}(S)\}$ or $n_{i}=2h$. Since $s>\operatorname{m}(S)$ and $h>\operatorname{m}(S)$ in both cases we obtain $n_{i}>\operatorname{m}(S)+h$, which is a contradiction. For the last statement, by Proposition~\ref{h,2n1} and the above discussion, if $\operatorname{m}(S)+h\geq n_{t-1}$ we have that $\operatorname{e}(T)\geq|\{n_{2},n_{3},\ldots,n_{t-1},2n_{1},h\}|=|\{n_{1},n_{2},\ldots,n_{t}\}|=\operatorname{e}(S)$.
\end{proof}

\begin{thm}
Let $S$ be a non special numerical semigroup and $T=\mathcal{A}(S)$. If there exists $x\in [\operatorname{F}(S)-\operatorname{m}(S)+1,\operatorname{F}(S)[\ \cap \operatorname{H}(S)$ then $\operatorname{e}(T)\geq\operatorname{e}(S)$. \label{I1}
\end{thm}
\begin{proof}
Let $\{\operatorname{m}(S)=n_{1}<n_{2}<\ldots<n_{t}\}$ be the set of minimal generators of $S$ and $I_1=[\operatorname{F}(S)-\operatorname{m}(S)+1,\operatorname{F}(S)[$. Let $x\in I_{1}\cap \operatorname{H}(S)$. First of all observe that if $s\in S$ then $x+s\geq x+n_{1}$. In particular $x+s\in I_{0}\subset S$. If $\operatorname{F}(S)<2\operatorname{m}(S)$ then $\operatorname{m}(S)\in I_{1}$ and since $S$ is not special, then we can choose $x=\max (\operatorname{SG}(S)\setminus \{\operatorname{F}(S)\})\in I_{1}$. If $\operatorname{F}(S)>2\operatorname{m}(S)$ then $n_{1}\in I_{j}$ with $j>1$, so $x>n_{1}$ and $2x>x+n_{1}\in I_{0}\subset S$, that is $x\in \operatorname{SG}(S)\setminus \{\operatorname{F}(S)\}$. Therefore, $h=\max (\operatorname{SG}(S)\setminus \{\operatorname{F}(S)\})\in I_{1}$. Consider the minimal generator $n_{t-1}$ of $S$. If $n_{t-1}\in I_{j}$ with $j\geq 1$ then $n_{1}+h>n_{t-1}$. If $n_{t-1}\in I_{0}$, let $k=n_{t-1}-n_{1}\in I_{1}$. In such a case, if $\operatorname{F}(S)>2\operatorname{m}(S)$ then $k\in \operatorname{SG}(S)\setminus \{\operatorname{F}(S)\}$, so $n_{1}+h\geq n_{1}+k=n_{t-1}$. If $\operatorname{F}(S)<2\operatorname{m}(S)$ we have that either $k<\operatorname{m}(S)$ or $k>\operatorname{m}(S)$, but in both cases $k\leq h$, since $h>\operatorname{m}(S)$ and $k$ is a special gap in the case $k>\operatorname{m}(S)$. In every case we obtain that $n_{1}+h\geq n_{t-1}$.
By Lemma~\ref{<nr} we have that $\operatorname{e}(T)\geq \operatorname{e}(S)$. 
\end{proof}

\noindent By the previous proposition to study the increasing or decreasing of the embedding dimension, with respect to $\mathcal{A}$, it remains to consider the semigroups with $[\operatorname{F}(S)-\operatorname{m}(S)+1,\operatorname{F}(S)[\subset S$. In such a case it may occur $\operatorname{e}(\mathcal{A}(S))<\operatorname{e}(S)$. Following the notation of Lemma~\ref{<nr}, these occurrences can be found when $n_{1}+h<n_{t-1}$. We give an example of this fact:

\begin{exa} \rm

Let $S$ be the numerical semigroup generated by the set $[761,768]\cup [11546,12305]$ 
We can do computation on such a semigroup using the package \texttt{numericalsgps} \cite{numericalsgps} in the computer algebra system \texttt{GAP}  \cite{GAP}.
\begin{lstlisting}
gap> G:=Concatenation([761..768],[11546..12305]);;
gap> s:=NumericalSemigroup(G);
gap> Length(SpecialGaps(s));
648
gap> h:=SpecialGaps(s)[647];
11537
gap> t:=AddSpecialGapOfNumericalSemigroup(11537,s);
<Numerical semigroup>
gap> t:=RemoveMinimalGeneratorFromNumericalSemigroup(761,t);
<Numerical semigroup with 762 generators>
gap> EmbeddingDimension(s);
655
gap> EmbeddingDimension(t);
652
\end{lstlisting}
So, in such a case, $\max (\operatorname{SG}(S)\setminus \{\operatorname{F}(S)\})=11537$. If $T=\mathcal{A}(S)$ then $652=\operatorname{e}(T)<\operatorname{e}(S)=655$. We observe also that $\operatorname{m}(S)+h=n_{t-7}$, where $\{n_{1}<n_{2}<\ldots <n_{t}\}$ is the set of minimal generators of $S$.
\label{examplei}
\end{exa}

We are able to provide a general pattern to obtain numerical semigroups having the same behaviour described in the above example, regarding the embedding dimension. This will be the content of a forthcoming paper.\\
With the next result we highlight a particular occurrence for a numerical semigroup $S$ for which $\operatorname{e}(T)>\operatorname{e}(S)$, with $T=\mathcal{A}(S)$.

\begin{prop} Let $S$ be a numerical semigroup, $S_{1}=\mathcal{A}(S)$, $S_{2}=\mathcal{A}(S_{1})=\mathcal{A}^{2}(S)$. Then $\operatorname{e}(S_{2})>\operatorname{e}(S_{1})$. \label{sporadic} \end{prop}

\begin{proof}
We consider that $\operatorname{F}(S)-\operatorname{m}(S)\in \operatorname{H}(S_{1})$, since $\operatorname{F}(S)-\operatorname{m}(S)\notin \operatorname{SG}(S)$. Moreover $\operatorname{m}(S_{1})>\operatorname{m}(S)$ and $\operatorname{F}(S_{1})=\operatorname{F}(S)$, so $\operatorname{F}(S)-\operatorname{m}(S)\in [\operatorname{F}(S_{1})-\operatorname{m}(S_{1})+1, \operatorname{F}(S_{1})[$. Then, if we consider $S_{2}=\mathcal{A}(S_{1})$, we have $\operatorname{e}(S_{2})>\operatorname{e}(S_{1})$ by Theorem~\ref{I1}.   
\end{proof}

\section{A family of trees related to the transfom $\mathcal{A}$}

In this section we show how to arrange all non special numerical semigroups with a given genus and a fixed number of left elements in a rooted tree, using the transform $\mathcal{A}$. This fact provides a tool to produce algorithmically these semigroups and to show a property linked to Wilf's conjecture.

\begin{defn} \rm
Let $g\in \mathbb{N}$ with $g>3$ and $n\in [2,g]$. We denote by $S_{g,n}$ the almost-ordinary numerical semigroup with genus $g$ and $\operatorname{n}(S_{g,n})=n$. Furthermore let $\mathcal{H}_{g,n}$ be the set whose elements are $S_{g,n}$ and all non special numerical semigroups $S$ such that $\operatorname{g}(S)=g$ and $\operatorname{n}(S)=n$.
\end{defn}

\begin{rmk}\rm
Let $S$ be a numerical semigroup of genus $g$. If $\operatorname{n}(S)=g$ then $\operatorname{F}(S)+1=2g$, so $S$ is symmetric. If $\operatorname{n}(S)=g-1$ then $\operatorname{F}(S)+2=2g$, so $S$ is pseudo-symmetric. In both cases $S$ is irreducible, in particular $\mathcal{H}_{g,g}=\{S_{g,g}\}$ and $\mathcal{H}_{g,g-1}=\{S_{g,g-1}\}$.
\end{rmk}


\begin{lemma}
Let $g\in \mathbb{N}$ with $g>3$ and $n\in [2,g-2]\cap \mathbb{N}$ and let $S\in \mathcal{H}_{g,n}$ with $S\neq S_{g,n}$. Consider the semigroup $T=(S\cup \{h\})\setminus \{\operatorname{m}(S)\}$ with $h=\max (\operatorname{SG}(S)\setminus \{\operatorname{F}(S)\})$. Then $\operatorname{m}(S)\in \operatorname{SG}(T)$ with $\operatorname{m}(S)<\operatorname{m}(T)$ and $h$ is a minimal generator of $T\cup \{\operatorname{m}(S)\}$  such that:
\begin{itemize}
\item If $S\cup \{h\}$ is irreducible then $\frac{\operatorname{F}(T)}{2}<h<\operatorname{F}(T)$
\item If $S\cup \{h\}$ is not irreducible then $\max(\operatorname{SG}(T\cup \{\operatorname{m}(S)\})\setminus \{\operatorname{F}(T)\})<h<\operatorname{F}(T)$. 
\end{itemize}
\label{treeNS1}
\end{lemma}
\begin{proof}
Since $h$ is a special gap of $S$, then $T\cup \{\operatorname{m}(S)\}=S\cup \{h\}$ is a numerical semigroup, that is $\operatorname{m}(S)\in \operatorname{SG}(T)$ and  $\operatorname{m}(S)<\operatorname{m}(T)$ by $h>\operatorname{m}(S)$. Trivially $h<\operatorname{F}(T)$ and $h>\frac{\operatorname{F}(S)}{2}=\frac{\operatorname{F}(T)}{2}$ by Lemma~\ref{max=max}, so if $S\cup \{h\}=T\cup \{\operatorname{m}(S)\}$ is irreducible we conclude. Otherwise let $y=\max(\operatorname{SG}(T\cup \{\operatorname{m}(S)\})\setminus \{\operatorname{F}(T)\})$, then $y\neq h$, $2y\neq \operatorname{F}(T\cup \{\operatorname{m}(S)\})=\operatorname{F}(S\cup \{h\})=\operatorname{F}(S)$ and $\operatorname{F}(T\cup \{\operatorname{m}(S)\})-y\notin T\cup \{\operatorname{m}(S)\}=S\cup \{h\}$, in particular $\operatorname{F}(S)-y\notin S$. This means that $y\in \{x\notin S\mid \operatorname{F}(S)-x\notin S, x\neq \frac{\operatorname{F}(S)}{2}\}$ but $h=\max\{x\notin S\mid \operatorname{F}(S)-x\notin S, x\neq \frac{\operatorname{F}(S)}{2}\}$ by Lemma~\ref{max=max}, so $h>y$.  
\end{proof}

\begin{lemma}
Let $T$ be a non ordinary numerical semigroup such that there exists $h\in \operatorname{SG}(T)$ with $h<\operatorname{m}(T)$ and let $y$ be a minimal generator of $T\cup \{h\}$ such that $y\neq h$ and:
\begin{itemize}
\item If $T\cup \{h\}$ is irreducible then $\frac{\operatorname{F}(T)}{2}<y<\operatorname{F}(T)$.
\item If $T\cup \{h\}$ is not irreducible then $\max(\operatorname{SG}(T\cup \{h\})\setminus \{\operatorname{F}(T)\})<y<\operatorname{F}(T)$. 
\end{itemize} 
Then $T=\mathcal{A}(S)$, with $S=(T\cup \{h\})\setminus \{y\}$. 
\label{treeNS2}
\end{lemma}

\begin{proof}
Let $S$ be the semigroup as above. We have $T=(S\cup \{y\})\setminus \{h\}$, so we have to prove that $h=\operatorname{m}(S)$ and $y=\max (\operatorname{SG}(S)\setminus \{\operatorname{F}(S)\})$. Since $h<\operatorname{m}(T)$ then $h<s$ for all $s\in T$. It follows that $h=\operatorname{m}(S)$. Since $T$ is not ordinary $h<\operatorname{F}(T)$, moreover $y<\operatorname{F}(T)$, so we have $\operatorname{F}(T)=\operatorname{F}(T\cup \{h\})=\operatorname{F}(S\cup \{y\})=\operatorname{F}(S)$. Put $X=\{ x\notin S\mid \operatorname{F}(S)-x\notin S, x\neq \frac{\operatorname{F}(S)}{2}\}$, it suffices, by Lemma~\ref{max=max}, to prove that $y=\max X  $. 
Observe that $T\cup \{h\}=S\cup \{y\}$ is a numerical semigroup and $y$ is a minimal generator of $T\cup \{h\}$, so $y\in \operatorname{SG}(S)\setminus \{\operatorname{F}(S)\}$ that is $y \in X$. Suppose $z=\max X$ and $z>y$. By Lemma~\ref{max=max}, $z\in \operatorname{SG}(S)$, $z\neq \operatorname{F}(S)=\operatorname{F}(T)=\operatorname{F}(T\cup\{h\})$ and $z>\frac{\operatorname{F}(S)}{2}=\frac{\operatorname{F}(T\cup \{h\})}{2}$. We prove that $z\in \operatorname{SG}(T\cup \{h\})$. Observe that for all $n\in \mathbb{N}$ then $z+n\neq y$. So let $x\in T\cup \{h\}$, we prove that $z+x\in T\cup\{h\}=S\cup \{y\}$. We have $z+h\in T\cup \{h\}$ since $h\in S$. Let $x\in T$, if $x\neq y$ then $x\in S$ and $z+x\in T\cup \{h\}$. Moreover $z+y\in T\cup \{h\}$ since by hypotheses $y>\max (\operatorname{SG}(T\cup \{h\})\setminus \{\operatorname{F}(T)\})>\frac{\operatorname{F}(T\cup \{h\})}{2}$ or $y>\frac{\operatorname{F}(T)}{2}=\frac{\operatorname{F}(T\cup \{h\})}{2}$, so $z+y>\operatorname{F}(T\cup \{h\})$. Furthermore $2z\in S$ and $2z\neq y$, so $2z \in T\cup \{h\}$. Therefore $z\in \operatorname{SG}(T\cup \{h\})\setminus \{\operatorname{F}(T)\}$, but this is a contradiction, because if $T\cup \{h\}$ is irreducible we have $\operatorname{SG}(T\cup \{h\})\setminus \{\operatorname{F}(T)\}= \emptyset$, and if $T\cup \{h\}$ is not irreducible we have $z<y$. So $y=\max X=\max (\operatorname{SG}(S)\setminus\{\operatorname{F}(S)\})$.
\end{proof}

\begin{defn} \rm
Let $g\in \mathbb{N}$ with $g>3$ and $n\in [2,g-2]$. We define the oriented graph $\mathcal{T}_{g,n}=(\mathcal{H}_{g,n},\mathcal{E})$ where $\mathcal{E}$ is the set of all pairs $(S,\mathcal{A}(S))$. If $(S,T)\in \mathcal{E}$ we say that $S$ is a \emph{child} of $T$. A numerical semigroup without children is called a \emph{leaf}.
\label{defTree}
\end{defn}

\begin{thm}
Let $g\in \mathbb{N}$ with $g>3$ and $n\in [2,g-2]$. The graph $\mathcal{T}_{g,n}$ is a rooted tree where the root is the almost-ordinary semigroup $S_{g,n}$. Moreover if $T\in \mathcal{H}_{g,n}$ then the children of $T$ are the semigroups $(T \cup \{h\})\setminus \{y\}$ where $h\in \operatorname{SG}(T)$ with $h<\operatorname{m}(T)$ and $y$ is a minimal generator of $T\cup \{h\}$ such that $y\neq h$ and one of the following holds:
\begin{itemize}
\item If $T\cup \{h\}$ is irreducible then $\frac{\operatorname{F}(T)}{2}<y<\operatorname{F}(T)$.
\item If $T\cup \{h\}$ is not irreducible then $\max(\operatorname{SG}(T\cup \{h\})\setminus \{\operatorname{F}(T)\})<y<\operatorname{F}(T)$. 
\end{itemize} 
\label{treeNS}
\end{thm}

\begin{proof}
Let $T\in \mathcal{H}_{g,n}$, we define the following sequence: 
\begin{itemize}
\item $T_{0}=T$.
\item $T_{i+1}=\left\lbrace\begin{array}{ll}
         \mathcal{A}(T_{i}) & \mbox{if $T_{i}\neq S_{g,n}$} \\
         S_{g,n} & \mbox{otherwise}
        \end{array}
        \right.$
\end{itemize}
\noindent In particular $T_{i}=\mathcal{A}^{i}(T)$ for all $i$. By Theorem~\ref{rootNS} there exists a nonnegative integer $k$ such that $T_{k}=\mathcal{A}^{k}(T)=S_{g,n}$. So the edges $(T_{0},T_{1}),(T_{1},T_{2}),\ldots,(T_{k-1},T_{k})$ provide a path from $T$ to $S_{g,n}$, hence $\mathcal{T}_{g,n}$ is a rooted tree whose root is $S_{g,n}$.\\ 
Let $T_{h,y}=(T \cup \{h\})\setminus \{y\}$ be a numerical semigroup as described above in the statement of the theorem. By Lemma~\ref{treeNS2} every pair $(T_{h,y},T)$ is an edge of $\mathcal{T}_{g,n}$ for every possible choice of $h$ and $y$. So every semigroup $T_{h,y}$ is a child of $T$ and these semigroups are exactly all the children of $T$ by Lemma~\ref{treeNS1}.
\end{proof}

Let $g\in \mathbb{N}$ with $g>3$ and $n\in [2,g-2]$. Starting with the special numerical semigroup $S_{g,n}$ it is possible to produce all non special numerical semigroups $S$ with genus $g$ and $\operatorname{n}(S)=n$ using the process described in the previous theorem. Lemma~\ref{treeNS1} and Lemma~\ref{treeNS2} ensure that all these semigroups are produced without redundancy. The process stops when we obtain all the leaves of the rooted tree $\mathcal{T}_{g,n}$. The leaves of this graph can be recognized as described below:

\begin{coro}
A numerical semigroup $T\in \mathcal{H}_{g,n}$ is a leaf in $\mathcal{T}_{g,n}$ if and only if for all $h\in \operatorname{SG}(T)$ one has $h>\operatorname{m}(T)$ or for every minimal generator $y$ of $T\cup \{h\}$ with $y\neq h$ and $y<\operatorname{F}(T)$ one has:
\begin{itemize}
\item If $T\cup \{h\}$ is irreducible then $y<\frac{\operatorname{F}(T)}{2}$.
\item If $T\cup \{h\}$ is not irreducible then $y<\max(\operatorname{SG}(T\cup \{h\})\setminus \{\operatorname{F}(T)\})$.
\end{itemize}
\label{leaves}
\end{coro}

\begin{exa} \rm
We build here $\mathcal{T}_{8,4}$. Let $S_{8,4}=\{0,8,9,10,12,\rightarrow\}$. All numerical semigroups involved in $\mathcal{T}_{8,4}$ have the same Frobenius number, that is $F=11$. As usual we denote by $h$ the greatest special gap different from $F$ for the semigroups involved. We start by computing the children of $S_{8,4}$. The special gaps smaller than multiplicity are $\{4,5,6,7\}$. So we consider:
\begin{enumerate}
\item $S_{8,4}\cup \{4\}=\{0,4,8,9,10,12,\rightarrow \}=\langle 4,9,10,15\rangle$, in which $h=6$ so by Theorem~\ref{treeNS} we have to consider the minimal generators in the set $\{9,10\}$.
\item $S_{8,4}\cup \{5\}=\{0,5,8,9,10,12,\rightarrow \}=\langle 5,8,9,12\rangle$, in which $h=7$ so by Theorem~\ref{treeNS} we have to consider the minimal generators in the set $\{8,9\}$.
\item $S_{8,4}\cup \{6\}=\{0,6,8,9,10,12,\rightarrow \}=\langle 6,8,9,10,13\rangle$, in which $h=7$ so by Theorem~\ref{treeNS} we have to consider the minimal generators in the set $\{8,9,10\}$.
\item $S_{8,4}\cup \{7\}=\{0,7,8,9,10,12,\rightarrow \}=\langle 7,8,9,10,12\rangle$, in which $h=6$ so by Theorem~\ref{treeNS} we have to consider the minimal generators in the set $\{8,9,10\}$.
\end{enumerate}  
The children of $S_{8,4}$ are the following:
\begin{itemize}
\item $S_{1}=(S_{8,4}\cup \{4\})\setminus \{9\}=\{0,4,8,10,12,\rightarrow\}$.
\item $S_{2}=(S_{8,4}\cup \{4\})\setminus \{10\}=\{0,4,8,9,12,\rightarrow\}$.
\item $S_{3}=(S_{8,4}\cup \{5\})\setminus \{8\}=\{0,5,9,10,12,\rightarrow\}$.
\item $S_{4}=(S_{8,4}\cup \{5\})\setminus \{9\}=\{0,5,8,10,12,\rightarrow\}$.
\item $S_{5}=(S_{8,4}\cup \{6\})\setminus \{8\}=\{0,6,9,10,12,\rightarrow\}$.
\item $S_{6}=(S_{8,4}\cup \{6\})\setminus \{9\}=\{0,6,8,10,12,\rightarrow\}$.
\item $S_{7}=(S_{8,4}\cup \{6\})\setminus \{10\}=\{0,6,8,9,12,\rightarrow\}$.
\item $S_{8}=(S_{8,4}\cup \{7\})\setminus \{8\}=\{0,7,9,10,12,\rightarrow\}$.
\item $S_{9}=(S_{8,4}\cup \{7\})\setminus \{9\}=\{0,7,8,10,12,\rightarrow\}$.
\item $S_{10}=(S_{8,4}\cup \{7\})\setminus \{10\}=\{0,7,8,9,12,\rightarrow\}$.
\end{itemize}
The numerical semigroups $S_{1},S_{2},S_{3},S_{4},S_{7}$ are leaves since their special gaps are smaller than their multiplicities. Also $S_{6}$ is a leaf, since $4$ is the unique special gap smaller than multiplicity and $S_{6}\cup \{4\}=\langle 4,6,13,15\rangle$ and it has $h=9$. By the same argument of $S_{6}$ also $S_{10}$ is a leaf.\\
By the same procedure we can check that $S_{5}$ has  the following child:
\begin{itemize}
\item $S_{11}=(S_{5}\cup \{3\})\setminus \{10\}=\{0,3,6,9,12,\rightarrow\}$.
\end{itemize}

The children of $S_{8}$ are the following:

\begin{itemize}
\item $S_{12}=(S_{8}\cup \{5\})\setminus \{9\}=\{0,5,7,10,12,\rightarrow\}$.
\item $S_{13}=(S_{8}\cup \{6\})\setminus \{9\}=\{0,6,7,10,12,\rightarrow\}$.
\item $S_{14}=(S_{8}\cup \{6\})\setminus \{10\}=\{0,6,7,9,12,\rightarrow\}$.
\end{itemize}

$S_{9}$ has one child:
\begin{itemize}
\item $S_{15}=(S_{9}\cup \{6\})\setminus \{10\}=\{0,6,7,8,12,\rightarrow\}$.
\end{itemize}

The numerical semigroups $S_{11},S_{12},S_{13},S_{14},S_{15}$ are leaves in $\mathcal{T}_{8,4}$. The graph $\mathcal{T}_{8,4}$ is pictured in Figure~\ref{FigureTree}.
\begin{center}
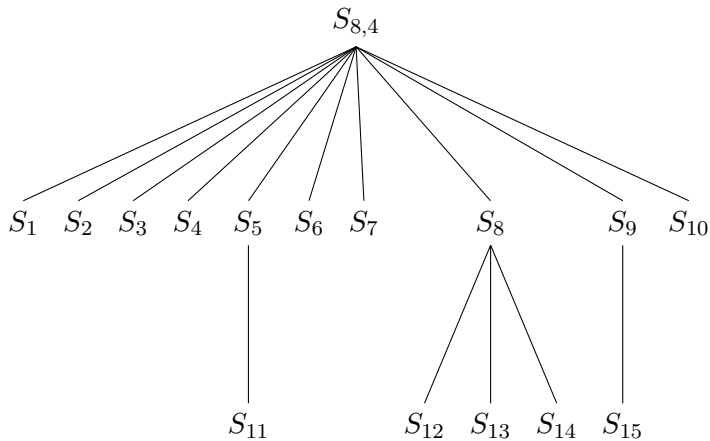
\begin{figure}
\begin{tikzpicture}
\tikzset{level distance=7em}
\Tree
	[.$S_{8,4}$  
		$S_1$ 
		$S_2$ 
		$S_3$
		$S_4$
		[.$S_5$ $S_{11}$ ] 
		$S_{6}$
		$S_{7}$
			[.$S_{8}$ $S_{12}$ $S_{13}$ $S_{14}$ ]
		[.$S_9$ $S_{15}$ ] 		
		$S_{10}$
	]
\end{tikzpicture}
\caption{The tree $\mathcal{T}_{8,4}$}
\label{FigureTree}
\end{figure}
\end{center}
\end{exa}

\vspace{3pt}
The particular behaviours of the transform $\mathcal{A}$ allow us to state the following consequences for Wilf's conjecture.

\begin{thm}
Let $\mathcal{F}_{g,n}$ be the set of all leaves of the tree $\mathcal{T}_{g,n}$. Suppose that:
\begin{enumerate}
\item All $S\in \mathcal{F}_{g,n}$ satisfy Wilf's conjecture.
\item $\mathcal{A}(S)$ satisfies Wilf's conjecture for all $S\in \mathcal{F}_{g,n}$ such that $[\operatorname{F}(S)-\operatorname{m}(S)+1,\operatorname{F}(S)[\subseteq S$
\end{enumerate}
Then all semigroups in $\mathcal{H}_{g,n}$ satisfy Wilf's conjecture.
\label{scheme}
\end{thm}

\begin{proof}
The transform $\mathcal{A}$ fixes all invariants involved in Wilf's conjecture except for the embedding dimension. Moreover, by Proposition~\ref{sporadic}, applying the transform $\mathcal{A}$, the embedding dimension is always increasing starting from a numerical semigroup $T$ such that $T=\mathcal{A}(S)$ with $S\in \mathcal{F}_{g,n}$. Furthermore if $S\in \mathcal{F}_{g,n}$ and there exists $h\in [\operatorname{F}(S)-\operatorname{m}(S)+1,\operatorname{F}(S)[ \cap \operatorname{H}(S)$, then from Theorem~\ref{I1}, $\mathcal{A}(S)$ has greater embedding dimension.
\end{proof}

\begin{prop}
Suppose there exists a numerical semigroup $S$ not satisfying Wilf's conjecture. Suppose $S\in \mathcal{H}_{g,n}$ and let $\mathcal{F}_{g,n}$ be the set of all leaves of the tree $\mathcal{T}_{g,n}$. Then there exists a numerical semigroup $T$ not satisfying Wilf's conjecture such that $T\in \mathcal{F}_{g,n}$ or $T=\mathcal{A}(T')$ for some $T'\in \mathcal{F}_{g,n}$.

\end{prop}

\begin{proof}
Easy consequence of the previous result.
\end{proof}

\begin{rmk} \rm
Observe that even if a numerical semigroup $S$ has a child in $\mathcal{T}_{g,n}$ with larger embedding dimension than $S$, it could be possible that $S$ has another child but with smaller embedding dimension than $S$. In particular we ask: does each $S\in \mathcal{H}_{g,n}$ have a child $T$ in $\mathcal{T}_{g,n}$ such that $\operatorname{e}(S)\geq \operatorname{e}(T)$?\\
\noindent If the previous question is true then Theorem~\ref{scheme} can be improved removing condition (2).
\label{improv1}
\end{rmk}

\begin{rmk} \rm
The trees $\mathcal{H}_{g,n}$ may be also used to generate all numerical semigroups of fixed genus, considering also an algorithm to generate all irreducible numerical semigroups of a fixed genus, that one can find for instance in \cite{blanco2013tree}. Anyway, we note that in \cite{blanco2012set} there exists a similar way to generate all numerical semigroups of a fixed genus, using different kinds of trees. Another way to achieve the same goal is described in the next Section.
\end{rmk}

\section{Another transform and the related trees} 

Let $S$ be a non ordinary numerical semigroup. The Frobenius number of the semigroup $S\cup \{\operatorname{F}(S)\}$ is defined as the \emph{sub-Frobenius number of $S$} and we denote it by $\operatorname{u}(S)$. In particular $$\operatorname{u}(S)=\max( \operatorname{H}(S)\setminus \{\operatorname{F}(S)\})$$ \\
In \cite{Bras-Amoros2018inproc-Different} the following transform is defined:

\begin{defn}\rm Let $\mathcal{Q}$ be the set of all non ordinary and non almost-ordinary numerical semigroups. $\mathcal{B}:\mathcal{Q}\rightarrow \mathcal{S}$ is a transform defined by:

$$\mathcal{B}(S)=(S\cup \{\operatorname{u}(S)\})\setminus \{\operatorname{m}(S)\}$$


\end{defn}

As mentioned in \cite{Bras-Amoros2018inproc-Different}, if $S\in \mathcal{Q}$ it is easy to show that there exists $n\in \mathbb{N}$ such that $\mathcal{B}^{n}(S)$ is almost-ordinary. 
The two transforms $\mathcal{A}$ and $\mathcal{B}$ are quite similar, observe for instance that also the transform $\mathcal{B}$ fixes all invariants involved in Wilf's conjecture except for the embedding dimension $\operatorname{e}(S)$. We pointed out some differences.

\begin{rmk} \rm
Let $S\in\mathcal{Q}$, observe that if $[\operatorname{F}(S)-\operatorname{m}(S)+1,\operatorname{F}(S)[\ \cap \operatorname{H}(S)\neq \emptyset$ then $\operatorname{u}(S)\in [\operatorname{F}(S)-\operatorname{m}(S)+1,\operatorname{F}(S)[$ and it is a special gap, since $\operatorname{u}(S)+\operatorname{m}(S)>\operatorname{F}(S)$. Moreover there are not gaps of $S$ different from $\operatorname{F}(S)$ and greater than $\operatorname{u}(S)$. In particular we have:
\begin{itemize}
\item[1)] $\operatorname{u}(S)=\max (\operatorname{SG}(S)\setminus \{\operatorname{F}(S)\})$  if $[\operatorname{F}(S)-\operatorname{m}(S)+1,\operatorname{F}(S)[\ \cap \operatorname{H}(S)\neq \emptyset$.
\item[2)] $\operatorname{u}(S)=\operatorname{F}(S)-\operatorname{m}(S)$ if $[\operatorname{F}(S)-\operatorname{m}(S)+1,\operatorname{F}(S)[\subset S$.
\end{itemize}
Therefore in the first case we have $\mathcal{B}(S)=\mathcal{A}(S)$. The second condition shows where the difference between $\mathcal{A}$ and $\mathcal{B}$ occurs.
\label{remBra}
\end{rmk}

\begin{exa}\rm
Let $S=\langle 6,11,13,15,16 \rangle= \mathbb{N}\setminus \{1,2,3,4,5,7,8,9,10,14,20\}$. In such a case:
\begin{itemize}
\item $\mathcal{A}(S)=(S\cup \{9\})\setminus \{6\}$
\item $\mathcal{B}(S)=(S\cup \{14\})\setminus \{6\}$
\end{itemize}
Observe that $S\cup \{14\}$ is not a numerical semigroup, since $6+14=20\notin S\cup \{14\}$. We note that $\mathcal{B}(S)$ is a numerical semigroup and that $\mathcal{B}(S)=\mathcal{A}^{2}(S)$. 

\label{exaBraDiff}
\end{exa}

If $S\in \mathcal{Q}$, from the previous example we observe that even if $(S\cup \{\operatorname{u}(S)\})\setminus \{\operatorname{m}(S)\}$ is a numerical semigroup it can occur that $S\cup \{\operatorname{u}(S)\}$ is not a semigroup. Actually the transform $\mathcal{B}$ can be defined equivalently as $\mathcal{B}(S)=(S\setminus \{\operatorname{m}(S)\})\cup \{\operatorname{u}(S)\}$, in particular $S\setminus \{\operatorname{m}(S)\}$ is a numerical semigroup.\\

As shown previously, once a transform is introduced for numerical semigroups then it is possible to arrange the
numerical semigroups into a graph. Now we describe how $\mathcal{B}$ allows one to arrange the set of numerical semigroups in a family of rooted trees following the procedure described in \cite{Bras-Amoros2018inproc-Different}.\\ 
We denote by $\mathcal{N}_{g,n}$ the set of all numerical semigroups with genus $g$ and number of left elements $n$.

\begin{defn} \rm
Let $g\in \mathbb{N}$ with $g>3$ and $n\in [2,g]$. We define the oriented graph $\mathcal{T}'_{g,n}=(\mathcal{N}_{g,n},\mathcal{E}')$ where $\mathcal{E}'$ is the set of all pairs $(S,\mathcal{B}(S))$.
\end{defn}

As in Definition~\ref{defTree}, if $(S,T)\in \mathcal{E}'$ we say that $S$ is a \emph{child} of $T$ and a numerical semigroup without children is called a \emph{leaf}. Recall that we denote $S_{g,n}$ the almost-ordinary numerical semigroup of genus $g$ and $\operatorname{n}(S_{g,n})=n$.

\begin{thm}[\cite{Bras-Amoros2018inproc-Different}]
Let $g\in \mathbb{N}$ with $g>3$ and $n\in [2,g]$. The graph $\mathcal{T}'_{g,n}$ is a rooted tree where the root is the almost-ordinary semigroup $S_{g,n}$. Moreover if $T\in \mathcal{N}_{g,n}$ then the children of $T$ in $\mathcal{T}'_{g,n}$ are the semigroups $(T \setminus \{y\})\cup \{h\}$ where $y$ is a minimal generator of $T$ such that $\operatorname{u}(T)<y<\operatorname{F}(T)$ and $h\in \operatorname{SG}(T\setminus \{y\})$ with $h<\operatorname{m}(T)$.
\label{treeNSBras}
\end{thm}

\begin{proof}
Similarly to Theorem~\ref{treeNS} one can prove that $\mathcal{T'}_{g,n}$ is a rooted tree whose root is $S_{g,n}$.\\ 
Let $T_{y,h}=(T\setminus \{y\})\cup \{h\}$ be a numerical semigroup as described above in the statement of the theorem. Then $\operatorname{u}(T_{y,h})=y$ and $\operatorname{m}(T_{y,h})=h$, that is $\mathcal{B}(T_{y,h})=T$, and $T_{y,h}$ is a child of $T$. Moreover, if $S$ is a child of a numerical semigroup $T$ then $T=(S\setminus \{\operatorname{m}(S)\})\cup \{\operatorname{u}(S)\}$, that is, $S= (T \setminus \{\operatorname{u}(S)\})\cup \{\operatorname{m}(S)\}$. Furthermore $\operatorname{u}(T)<\operatorname{u}(S)<\operatorname{F}(T)$ and $\operatorname{m}(S)\in \operatorname{SG}(T\setminus \{\operatorname{u}(S)\})$, in particular $\operatorname{m}(S)<\operatorname{m}(T)$.
\end{proof}

\begin{coro}
A numerical semigroup $T\in \mathcal{N}_{g,n}$ is a leaf in $\mathcal{T}'_{g,n}$ if and only if  $T$ has no minimal generators in the interval $[\operatorname{u}(S),\operatorname{F}(S)]$ or for all minimal generators $y\in [\operatorname{u}(S),\operatorname{F}(S)]$ and for all $h\in \operatorname{SG}(T\setminus \{y\})$ one has $h>\operatorname{m}(T)$.
\label{leavesBra}
\end{coro}

Unlike $\mathcal{A}$, the transform $\mathcal{B}$ has monotone behaviour with respect to the growth of the embedding dimension.

\begin{prop}
Let $S$ be a non ordinary and non almost-ordinary numerical semigroup. Then $\operatorname{e}(\mathcal{B}(S))\geq \operatorname{e}(S)$
\label{I2}
\end{prop}
\begin{proof}
If $[\operatorname{F}(S)-\operatorname{m}(S)+1,\operatorname{F}(S)[\ \cap \operatorname{H}(S)\neq \emptyset$ then $\mathcal{B}(S)=\mathcal{A}(S)$ and we obtain the claim by Theorem~\ref{I1}. So we can assume that $[\operatorname{F}(S)-\operatorname{m}(S)+1,\operatorname{F}(S)[ \subset S$, in particular $\mathcal{B}(S)= (S\cup \{\operatorname{F}(S)-\operatorname{m}(S)\})\setminus \{\operatorname{m}(S)\}$. Observe that $2\operatorname{m}(S)$ and $\operatorname{F}(S)-\operatorname{m}(S)$ are minimal generators of $\mathcal{B}(S)$, since $S$ is not almost ordinary. Let $\{\operatorname{m}(S)=n_1<n_2<\ldots<n_t\}$ be the set of minimal generators of $S$, in this case we have $n_{t-1}<\operatorname{F}(S)$  otherwise $n_{t-1}-\operatorname{m}(S)\in S$. Suppose that $n_i$ is not a minimal generator of $\mathcal{B}(S)$ for some $i\in\{1,\ldots,t-1\}$, then $n_{i}=(\operatorname{F}(S)-\operatorname{m}(S))+t$ for some $t\in S\setminus \{0\}$, in particular $n_i>\operatorname{F}(S)$, which is a contradiction. So $\operatorname{e}(\mathcal{B}(S))\geq |\{2\operatorname{m}(S), \operatorname{F}(S)-\operatorname{m}(S), n_2,\ldots, n_{t-1}\}|=|\{n_1,n_2,\ldots,n_t\}|=\operatorname{e}(S)$.
\end{proof}

\begin{coro}
Let $\mathcal{F}'_{g,n}$ be the set of all leaves of the tree $\mathcal{T}'_{g,n}$ and suppose that all $S\in \mathcal{F}'_{g,n}$ satisfy Wilf's conjecture. Then all semigroups in $\mathcal{N}_{g,n}$ satisfy Wilf's conjecture.
\label{schemeBras}
\end{coro}

\begin{prop}
Suppose there exists a numerical semigroup $S$ not satisfying Wilf's conjecture. Suppose $S\in \mathcal{N}_{g,n}$ and let $\mathcal{F}'_{g,n}$ be the set of all leaves of the tree $\mathcal{T}'_{g,n}$. Then there exists a numerical semigroup in $\mathcal{F}'_{g,n}$  not satisfying Wilf's conjecture.

\end{prop}


Notice that applying $\mathcal{B}$ the embedding dimension is always increasing, so a counterexample to Wilf's conjecture will lead to a counterexample among the leaves of a tree $\mathcal{T}'_{g,n}$. Theorem~\ref{treeNSBras} allows one to compute the trees $\mathcal{T}'_{g,n}$ algorithmically. Moreover, all numerical semigroups of a fixed genus can be computed by this procedure and, unlike $\mathcal{A}$, also the set of all irreducible numerical semigroups of a fixed genus $g$ can be computed since they are the numerical semigroups in the trees $\mathcal{T}'_{g,g-1}$ and $\mathcal{T}'_{g,g}$.

\begin{rmk}\rm 
It could be interesting to have a more complete comparison between the two transforms and between the trees $\mathcal{T}_{g,n}$ and $\mathcal{T}'_{g,n}$.\\
In particular we can ask if the set $\mathcal{F}_{g,n}$ of all leaves of $\mathcal{T}_{g,n}$ is comparable with the set $\mathcal{F}'_{g,n}$ of all leaves of $\mathcal{T}'_{g,n}$. That is, if it occurs $\mathcal{F}_{g,n}\subseteq \mathcal{F}'_{g,n}$ or $\mathcal{F}'_{g,n}\subseteq \mathcal{F}_{g,n}$ or neither for some $g$ and $n$.

\end{rmk}

\section{Concluding remarks and possible developments}

In this paper we have provided some results concerning Wilf's conjecture, that is, we reduce the study of Wilf's conjecture to the study of the set $\mathcal{F}'_{g,n}$ of all leaves of the tree $\mathcal{T}'_{g,n}$, with $12<n<\frac{g}{3}$, or considering the set $\mathcal{F}_{g,n}$ of all leaves of the tree $\mathcal{T}_{g,n}$ and to the transformed of some of them using $\mathcal{A}$. In fact, it is known by \cite[Proposition 3.15]{Delgado2020} that if $S$ is a numerical semigroup such that $\operatorname{n}(S)\leq 12$ then $S$ satisfies Wilf's conjecture, and the same occurs if $3\operatorname{n}(S)\geq \operatorname{g}(S)$ as a consequence of \cite[Corollary 2.7]{dobbs2003question}. This may be an interesting result from a theoretical and a computational point of view. Moreover, with computational methods, in \cite{fromentin2016exploring} it has been proved that every numerical semigroup of genus $g$ with $g\leq 60$ satisfies Wilf's conjecture. This bound has been improved up to $g=80$, as cited in \cite{eliahou2019wilf}.
Considering such results, if one wants to find an example of a semigroup not satisfying Wilf's conjecture, one may investigate the trees $\mathcal{T}_{g,n}$ and $\mathcal{T'}_{g,n}$ with $g> 80$ and with $12<n<\frac{g}{3}$, in particular not the whole trees but their set of leaves. The above results allow to compute these trees algorithmically. A possible development can be to study the trees $\mathcal{T}_{g,n}$ and $\mathcal{T}'_{g,n}$ from a computational point of view.\\
An important reduction was introduced in \cite{eliahou2017wilf}, where an invariant $\operatorname{E}(S)$ is associated to each numerical semigroup, studied also in \cite{delgado2018question}, where it is called an \emph{Eliahou number}. It is defined as $\operatorname{E}(S)=|Q|\operatorname{n}(S)-q|D|+\rho$, where $Q$ is the set of minimal generators smaller than $\operatorname{F}(S)$, $D$ is the set of the elements $x \in [\operatorname{F}(S)+1,\operatorname{F}(S)+\operatorname{m}(S)]$ such that $x$ is not a minimal generator, and $q, \rho\in \mathbb{N}$ are the only numbers such that  $\operatorname{F}(S)+1=q\operatorname{m}(S)-\rho$. It has been proved, in \cite{eliahou2017wilf}, that if $\operatorname{E}(S)\geq 0$ then $S$ satisfies Wilf's conjecture. So the study of Wilf's conjecture is reduced to those numerical semigroups with $\operatorname{E}(S)<0$. Some families with an infinite number of numerical semigroups satisfying $\operatorname{E}(S)<0$ are provided in \cite{delgado2018question} and \cite{eliahou2019near}. By \cite{eliahou2017wilf} we know also that the proportion of numerical semigroups with $\operatorname{E}(S)<0$, with respect to all numerical semigroups, is asymptotically small. Just to consider an interesting feature, the numerical semigroup with the smallest genus satisfying $\operatorname{E}(S)<0$ has genus $43$. So, with confidence, we can state that Eliahou's result is a much finer reduction for studying Wilf's conjecture. In fact, even if in this paper it is not investigated the asymptotic proportion of the semigroups in $\mathcal{F}_{g,n}$ and $\mathcal{F}'_{g,n}$, there is a huge number of them, since there exist such numerical semigroups for each value of the genus and it is possible that the trees $\mathcal{T}_{g,n}$ and $\mathcal{T}'_{g,n}$ may contain too many leaves compared with the non-leaves. We pose the following problem: where are the numerical semigroups $S$, such that $\operatorname{E}(S)<0$, located in $\mathcal{T}_{g,n}$ or $\mathcal{T}'_{g,n}$? Are there any of them (or all) in the set $\mathcal{F}'_{g,n}$ of the leaves of $\mathcal{T}'_{g,n}$? Or in the set $\mathcal{F}_{g,n}$ of the leaves of $\mathcal{T}_{g,n}$ or in their transformed by $\mathcal{A}$? \\

We also mention that in \cite{cistodelgadosanchez} the ordinarization transform (defined in the preliminaries of this paper) was introduced in a more general setting than numerical semigroups, that is for submonoids of $\mathbb{N}^{d}$ with finite complement in $\mathbb{N}^{d}$, a class of monoids introduced in \cite{failla2016algorithms}. For such a class of submonoids of $\mathbb{N}^d$, also a generalization of Wilf's conjecture was proposed in \cite{preprint} and in \cite{garcia2018extension}. So, it could be interesting to investigate if it is possible to formulate a generalization of the definitions of $\mathcal{A}$ and $\mathcal{B}$ for submonoids in $\mathbb{N}^{d}$.

\subsection*{Acknowledgements}
I wish to thank Manuel Delgado for bringing article \cite{Bras-Amoros2018inproc-Different} to my attention as well as for his comments and helpful suggestions.

\bibliographystyle{plain}
\bibliography{miabiblio}

\end{document}